\documentclass[letterpaper]{article}
\usepackage[margin = 1in]{geometry}
\usepackage{graphicx,color}
\usepackage[cp850]{inputenc}
\usepackage[T1]{fontenc}
\usepackage[english]{babel}
\usepackage{cite}
\usepackage{rotating}
\usepackage[f]{esvect}
\usepackage{amsmath, amsthm, amssymb}
\usepackage{rotating}
\usepackage{lscape}
\usepackage{mathrsfs}
\usepackage{mathdots}
\usepackage{tikz}

\def\tablenotes{\bgroup\parfillskip=0pt plus 1fil
\leftskip=0pt\relax \rightskip=0pt
\vskip2pt\footnotesize}
\def\endtablenotes{\vskip1pt\egroup}

\newtheorem{theorem}{Theorem}[section]

\newtheorem{lemma}[theorem]{Lemma}

\newtheorem{definition}[theorem]{Definition}
\newcommand{\sech}{\mathrm{sech}}

\newcommand{\diag}{\mathrm{diag} \,}

\allowdisplaybreaks

\begin{document}
\title{Centrosymmetric Matrices in the Sinc Collocation Method for Sturm-Liouville Problems}

\author{Philippe Gaudreau and Hassan Safouhi\footnote{Corresponding author: hsafouhi@ualberta.ca}~\footnote{This work is supported by the Natural Sciences and Engineering Research Council of Canada~(NSERC) - Grant 250223-2011.}\\
Mathematical Section, Campus Saint-Jean\\
University of Alberta\\
8406, 91 Street, Edmonton, Alberta T6C 4G9, Canada}

\date{}

\maketitle

\vspace*{1.25cm}

\vspace*{1.25cm}
{\bf \large Abstract}. \hskip 0.10cm

Recently, we used the Sinc collocation method with the double exponential transformation to compute eigenvalues for singular Sturm-Liouville problems. In this work, we show that the computation complexity of the eigenvalues of such a differential eigenvalue problem can be considerably reduced when its operator commutes with the parity operator. In this case, the matrices resulting from the Sinc collocation method are centrosymmetric. Utilizing well known properties of centrosymmetric matrices, we transform the problem of solving one large eigensystem into solving two smaller eigensystems. We show that only $\frac{1}{N+1}$ of all components need to be computed and stored in order to obtain all eigenvalues, where $2N+1$ corresponds to the dimension of the eigensystem.  We applied our result to the Schr\"odinger equation with the anharmonic potential and the numerical results section clearly illustrates the substantial gain in efficiency and accuracy when using the proposed algorithm.

\vspace*{1.25cm}
{\bf Keywords}

Sturm-Liouville eigenvalue problem. Schr\"odinger equation. Anharmonic oscillators. Double exponential Sinc-Collocation method. Centrosymmetry.

\clearpage
\section{Introduction}
In science and engineering, differential eigenvalue problems occur abundantly. Differential eigenvalue problems can arise when partial differential equations are solved using the method of separation of variables. Consequently, they also play an important role in Sturm-Liouville (SL) differential eigenvalue problems \cite{Amrein2005}. For example, the solution of the wave equation can be expressed as the sum of standing waves. The frequencies of these standing waves are precisely the eigenvalues of its corresponding Sturm-Liouville problem. Similarly, in quantum mechanics, the energy eigenvalues associated with a Hamiltonian operator are modelled using the time-independent Schr\"odinger equation which is in fact a special case of a Sturm-Liouville differential eigenvalue problem.

Recently, collocation and spectral methods have shown great promise for solving singular Sturm-Liouville differential eigenvalue problems~\cite{Auzinger2006,Chanane2007}. More specifically, the Sinc collocation method (SCM) \cite{Tharwat2013a,Tharwat2013,Jarratt1990a} has been shown to yield exponential convergence. During the last three decades the SCM has been used extensively to solve many problems in numerical analysis. The applications include numerical integration, linear and non-linear ordinary differential equations, partial differential equations, interpolation and approximations to functions~\cite{Stenger1981,Stenger2000}. The SCM applied to Sturm-Liouville problems consists of expanding the solution of a SL problem using a basis of Sinc functions. By evaluating the resulting approximation at the Sinc collocation points separated by a fixed mesh size $h$, one obtains a matrix eigenvalue problem or generalized matrix eigenvalue problem for which the eigenvalues are approximations to the eigenvalues of the SL operator. In~\cite{Safouhi50-arXiv}, we used the double exponential Sinc collocation method (DESCM) to compute the eigenvalues of singular Sturm-Liouville boundary value problems. The DESCM leads to a generalized eigenvalue problem where the matrices are symmetric and positive-definite. In addition, we demonstrate that the convergence of the DESCM is of the rate ${\cal O} \left( \frac{N^{5/2}}{\log(N)^{2}}  e^{-\kappa N/\log(N)}\right)$ for some $\kappa>0$ as $N\to \infty$, where $2N+1$ is the dimension of the resulting generalized eigenvalue system. The DESCM was also applied successfully to the Schr\"odinger equation with the anharmonic oscillators~\cite{Safouhi49}.

In the present contribution, we show how a parity symmetry of the Sturm-Liouville operator can be conserved and exploited when converting our differential eigenvalue problem into a matrix eigenvalue problem. Indeed, given certain parity assumptions, the matrices resulting from the DESINC method are not only symmetric and positive definite; they are also centrosymmetric. The study of centrosymmetry has a long history \cite{Saibel1942, Cantoni1976a, Cruse1977, A1984, Stuart1988, Hill1990, Muthiyalu1992, Nield1994}. However, the last two decades has stemmed much research focused on the properties and applications of centrosymmetric matrices ranging from iterative methods for solving linear equations to least-squares problems to inverse eigenvalue problems~\cite{Melman2000,Tao2002,Lu2002,Abu-jeib2002,Zhou2003a,Zhou2003,Liu2003,Fassbender2003,Trench2004,Trench2004a, Zhongyun2005,Liu2005,Tian2007,Li2012,El-Mikkawy2013}.

Using the eigenspectrum properties of symmetric centrosymmetric matrices presented in \cite{Cantoni1976a}, we apply the DESCM algorithm to Sturm-Liouville eigenvalue problems and demonstrate that solving the resulting generalized eigensystem of dimension $(2N+1)\times(2N+1)$ is equivalent to solving the two smaller eigensystems of dimension $N \times N$ and $(N+1) \times (N+1)$. Moreover, we also demonstrate that only $\frac{1}{N+1}$ of all components need to be stored at every iteration in order to obtain all generalized eigenvalues. To illustrate the gain in efficiency obtained by this method, we apply the DESCM method to the time independent Schr\"odinger equation with an anharmonic potential. Furthermore, it is worth mentioning that research concerning inverse eigenvalue problems where the matrices are assumed centrosymmetric has been the subject of much research recently \cite{Trench2004a, Zhou2003a}. Consequently, the combination of these results and our findings could lead to a general approach for solving inverse Sturm-Liouville problems.

All calculations are performed using the programming language Julia and all the codes are available upon request.

\section{Definitions and basic properties}~\label{GenDef}
The sinc function valid for all $z \in \mathbb{C}$ is defined by the following expression:
\begin{equation} \label{formula: sinc functions}
\textrm{sinc}(z) = \left\{ \begin{array}{cc} \dfrac{\sin(\pi z)}{\pi z} &\quad \textrm{for} \quad z \neq 0 \\[0.3cm]
1  &\quad \textrm{for} \quad z=0. \end{array} \right.
\end{equation}

For $j \in \mathbb{Z}$ and $h$ a positive number, we define the Sinc function $S(j,h)(x)$ by:
\begin{equation}\label{formula: Sinc function}
S(j,h)(x) = \textrm{sinc}\left( \dfrac{x-jh}{h}\right) \quad \textrm{for} \quad x \in \mathbb{C}.
\end{equation}

The Sinc function defined in \eqref{formula: Sinc function} form an interpolatory set of functions with the discrete orthogonality property:
\begin{equation}
S(j,h)(kh)  =  \delta_{j,k}   \qquad \textrm{for} \qquad j,k \in \mathbb{Z},
\end{equation}
where $\delta_{j,k}$ is the Kronecker delta function.

\begin{definition}\cite{Stenger1981}
Given any function $v$ defined everywhere on the real line and any $h>0$, the symmetric truncated Sinc expansion of $v$ is defined by the following series:
\begin{equation}\label{formula: y in sinc function}
C_{N}(v,h)(x)  = \displaystyle \sum_{j=-N}^{N} v_{j,h} \, S(j,h)(x),
\end{equation}
where $v_{j,h} = v(jh)$.
\end{definition}
The Sturm-Liouville (SL) equation in Liouville form is defined as follows:
\begin{align} \label{formula: sturm-liouville problem}
Lu(x) & =  - u^{\prime \prime}(x) + q(x) u(x)  \,=\,  \lambda \rho(x)u(x) \nonumber \\
&  \hskip -0.5cm a < x < b  \qquad \qquad  u(a) = u(b)=0,
\end{align}
where $ -\infty \leq a < b \leq \infty$. Moreover, we assume that the function $q(x)$ is non-negative and the weight function $\rho(x)$ is positive. The values $\lambda$ are known as the eigenvalues of the SL equation.

In \cite{Gaudreau2014a}, we apply the DESCM to obtain an approximation to the eigenvalues $\lambda$ of equation \eqref{formula: sturm-liouville problem}. We initially applied Eggert et al.'s transformation to equation \eqref{formula: sturm-liouville problem} since it was shown that the proposed change of variable results in a symmetric discretized system when using the Sinc collocation method \cite{Eggert1987}. The proposed change of variable is of the form~\cite[Defintion 2.1]{Eggert1987}:
\begin{equation}
v(x) = \left(\sqrt{ (\phi^{-1})^{\prime} } \, u \right) \circ \phi(x) \qquad \Longrightarrow \qquad  u(x)  =  \dfrac{ v \circ \phi^{-1}(x)}{\sqrt{ (\phi^{-1}(x))^{\prime}}},
\label{formula: EggertSub}
\end{equation}
where $\phi^{-1}(x)$ is a conformal map of a simply connected domain in the complex plane with boundary points $a\neq b$ such that $\phi^{-1}(a)=-\infty$ and $\phi^{-1}(b)=\infty$.

Applying the change of variable \eqref{formula: EggertSub} into equation~\eqref{formula: sturm-liouville problem}, one obtains~\cite{Eggert1987}:
\begin{equation}\label{formula: transformed sturm-liouville problem}
\mathcal{L} \, v(x) =  - v^{\prime \prime}(x) + \tilde{q}(x) v(x) = \lambda \rho(\phi(x))(\phi^{\prime}(x))^{2} v(x)\quad \textrm{with} \quad \lim_{|x| \to \infty} v(x) = 0,
\end{equation}
where:
\begin{equation}\label{formula: q tilde}
\tilde{q}(x)  =  - \sqrt{\phi^{\prime}(x)} \, \dfrac{{\rm d}}{{\rm d} x} \left( \dfrac{1}{\phi^{\prime}(x)} \dfrac{{\rm d}}{{\rm d} x}( \sqrt{\phi^{\prime}(x)})  \right) + (\phi^{\prime}(x))^{2} q(\phi(x)).
\end{equation}

To implement the double exponential transformation, we use a conformal mapping $\phi(x)$ such that the solution to equation \eqref{formula: transformed sturm-liouville problem} decays double exponentially. In other words, we need to find a function $\phi(x)$ such that:
\begin{equation}
|v(x)| \leq A \exp \left( -B \exp \left( \gamma |x| \right) \right),
\end{equation}
for some positive constants $A,B,\gamma$. Examples of such mappings are given in~\cite{Gaudreau2014a, Mori2001}.

Applying the SCM method, we obtain the following generalized eigenvalue problem:
\begin{align}\label{formula: matrix solution}
\mathcal{L} \, {\bf C}_{M}(v,h) &  = {\bf A}{\bf v} \,=\,  \mu {\bf D}^{2}{\bf v} \quad \Longrightarrow \quad ({\bf A} - \mu {\bf D}^{2} ){\bf v} \,=\, 0,
\end{align}
where the vectors ${\bf v}$ and ${\bf C}_{M}(v,h)$ are given by:
\begin{equation}
{\bf v} = (v_{-N,h},\ldots, v_{N,h})^{T} \qquad \textrm{and} \qquad {\bf C}_{N}(v,h)  = (C_{N}(v,h)(-Nh), \ldots, C_{N}(v,h)(Nh) )^{T},
\label{EQVECTORCN001}
\end{equation}
and $\mu$ are approximations of the eigenvalues $\lambda$ of equation \eqref{formula: transformed sturm-liouville problem}. For more details on the application of the SCM, we refer the readers to \cite{Safouhi50-arXiv}.

As in \cite{Stenger1997a}, we let $\delta^{(l)}_{j,k}$ denote the $l^{\rm th}$ Sinc differentiation matrix with unit mesh size:
\begin{equation}\label{formula: delta matrices}
\delta^{(l)}_{j,k} =  h^{l} \left. \left( \dfrac{\rm d}{{\rm d}x} \right)^{l} S(j,h)(x) \right|_{x=kh}.
\end{equation}

The entries $A_{j,k}$ of the $(2N+1) \times (2N+1)$ matrix ${\bf A}$ are then given by:
\begin{equation}\label{formula: A components}
A_{j,k} =   -\dfrac{1}{h^{2}} \, \delta^{(2)}_{j,k} + \tilde{q}(kh) \, \delta^{(0)}_{j,k} \qquad {\rm with} \qquad -N \leq j,k \leq N,
\end{equation}
and the entries $D^{2}_{j,k}$ of the $(2N+1) \times (2N+1)$ diagonal matrix ${\bf D}^{2}$ are given~by:
\begin{equation}\label{formula: D components}
D^{2}_{j,k} =  (\phi^{\prime}(kh))^{2} \rho(\phi(kh)) \, \delta^{(0)}_{j,k}  \qquad {\rm with} \qquad -N \leq j,k \leq N.
\end{equation}

As previously mentioned, Eggert et al.'s transformation leads to the matrices $\,{\bf A}$ and ${\bf D}^2$ to be symmetric and positive definite. However, as will be illustrated in the next section, given certain parity assumptions, these matrices yield even more symmetry.

\section{Centrosymmetric properties of the matrices ${\bf A}$ and ${\bf D}^2$}
In this section, we present some properties of the matrix $\,{\bf A}$ and ${\bf D}^2$ that will be beneficial in the computation of their eigenvalues. The matrices ${\bf A}$ and ${\bf D}^2$ are symmetric positive definite matrices when equation \eqref{formula: transformed sturm-liouville problem} is discretized using the Sinc collocation method. Additionally, given certain parity assumptions on the functions $q(x)$, $\phi(x)$ and $\rho(x)$ in equation \eqref{formula: transformed sturm-liouville problem}, the matrices ${\bf A}$ and ${\bf D}^2$ will also be centrosymmetric.

\begin{definition}\cite[Section 5.10]{Bransden2000}
Let $\mathcal{J}$ denote the parity operator defined by:
\begin{equation}
\mathcal{J}f(x) = f(-x),
\end{equation}
where $f(x)$ is a well defined function being acted upon by $\mathcal{J}$.
\end{definition}

\begin{definition}
An operator $\mathcal{B}$ is said to commute with parity operator $\mathcal{J}$ if it satisfies the following relation:
\begin{equation}
\mathcal{B}\mathcal{J}f(x) = \mathcal{J}\mathcal{B}f(x).
\end{equation}
Equivalently, we can say that the the commutator between $\mathcal{B}$ and $\mathcal{J}$ is zero, that is:
\begin{equation}
[\mathcal{B},\mathcal{J}] = \mathcal{B}\mathcal{J} - \mathcal{J}\mathcal{B} = 0.
\end{equation}
\end{definition}

\begin{definition}\cite[Definition 5]{Weaver1985}
An exchange matrix denoted by ${\bf J}$ is a square matrix with ones along the anti-diagonal and zeros everywhere else:
\begin{equation}
{\bf J} =
  \begin{pmatrix}
     \multicolumn{1}{c}{\emph{\text{\kern 0em\smash{\raisebox{-0.5ex}{\Large 0}}}}}&   & 1 \\
     & \iddots &  \\
    1 & & \multicolumn{1}{c}{\emph{\text{\kern0.5em\smash{\raisebox{0.5ex}{\Large 0}}}}}
  \end{pmatrix}.
\end{equation}
\end{definition}

\begin{definition}\cite[Definition 2]{Weaver1985}
Let {\bf B} be a matrix of dimension $(2N+1) \times (2N+1)$ with components $B_{j,k}$ for $-N \leq j,k \leq N$. ${\bf B}$ is centrosymmetric if and only if ${\bf B}$ satisfies the following property:
\begin{equation}\label{formula: centrosymmetric B}
{\bf B}{\bf J}={\bf J}{\bf B},
\end{equation}
where ${\bf J}$ is an exchange matrix of dimension $(2N+1) \times (2N+1)$.
Writing equation~\eqref{formula: centrosymmetric B} in a component form, we have the following relation:
\begin{equation}\label{formula: centrosymmetric B component}
B_{-j,-k} = B_{j,k} \quad \textrm{for} \quad -N\leq j,k\leq N.
\end{equation}
\end{definition}

We now present the following Theorem establishing the connection between symmetries of the Sturm-Liouville operator and its resulting matrix approximation.

\begin{theorem}\label{theorem: H centrosymmetric}
Let $\mathcal{L}$ denote the operator of the transformed Sturm-Liouville problem in equation \eqref{formula: transformed sturm-liouville problem}:
\begin{equation}\label{formula: transformed sturm-liouville problem operator}
\mathcal{L} = \dfrac{1}{\rho(\phi(x))(\phi^{\prime}(x))^{2}} \left( - \dfrac{ d^{2}}{dx^2} + \tilde{q}(x) \right).
\end{equation}
If the commutator $[\mathcal{L}, \mathcal{J}] = 0$, where $\mathcal{J}$ is the parity operator, then the matrices $\,{\bf A}$ and ${\bf D}^2$  defined by equations \eqref{formula: A components} and \eqref{formula: D components} resulting from the DESCM are centrosymmetric.
\end{theorem}

\underline{\bf Proof} \hskip 0.25cm The commutator $[\mathcal{L}, \mathcal{J}] = 0$ if and only if $q(x)$ and $\rho(x)$ are even functions and $\phi(x)$ is an odd function.

If $\phi(x)$ is an odd function, then $\phi^{\prime}(x)$ is even, $ \phi^{\prime \prime}(x)$ is odd and $ \phi^{\prime \prime \prime}(x)$ is even. From this and equation \eqref{formula: q tilde}, it follows that $\tilde{q}(x)$ is even.

In order to show that the resulting matrices $\,{\bf A}$ and ${\bf D}^2$ are centrosymmetric, we demonstrate that both these matrices satisfy equation \eqref{formula: centrosymmetric B component}. Before doing so, it is important to notice that the $l^{\rm th}$ Sinc differentiation matrices defined in equation \eqref{formula: delta matrices} have the following symmetric properties:
\begin{align}
\delta^{(l)}_{-j,-k} & =  h^{l} \left. \left( \dfrac{\rm d}{{\rm d}x} \right)^{l} S(-j,h)(x) \right|_{x=-kh}  \,= \, \begin{cases}
~~\delta^{(l)}_{j,k} & \textrm{if} \quad \textrm{$l$ is even} \\[0.1cm]
-\delta^{(l)}_{j,k} & \textrm{if} \quad \textrm{$l$ is odd}.
\end{cases}
\end{align}

Hence, the $l^{\rm th}$ Sinc differentiation matrices are centrosymmetric if $l$ is even. It is worth noting that when $l$ is odd, the Sinc differentiation matrices are skew-centrosymmetric \cite{Trench2004}. Consequently, investigating the form for the components of the matrix ${\bf A}$ in equation \eqref{formula: A components}, we obtain:
\begin{eqnarray}
A_{-j,-k} & = &   -\dfrac{1}{h^{2}} \, \delta^{(2)}_{-j,-k} + \tilde{q}(-kh) \, \delta^{(0)}_{-j,-k} \nonumber \\
          & = &   -\dfrac{1}{h^{2}} \, \delta^{(2)}_{j,k} + \tilde{q}(kh) \, \delta^{(0)}_{j,k} \nonumber\\
          & = & A_{j,k}.
\end{eqnarray}

Similarly, investigating the form for the components of the matrix ${\bf D}^2$ in equation \eqref{formula: D components}, we obtain:
\begin{eqnarray}
D^{2}_{-j,-k} & = &  (\phi^{\prime}(-kh))^{2} \rho(\phi(-kh)) \, \delta^{(0)}_{-j,-k}  \nonumber \\
              & = &  (\phi^{\prime}(kh))^{2} \rho(\phi(kh)) \, \delta^{(0)}_{j,k}  \nonumber \\
              & = &  D^{2}_{j,k}.
\end{eqnarray}

Both matrices ${\bf A}$ and ${\bf D}^2$ satisfy equation \eqref{formula: centrosymmetric B component}. From this it follows that ${\bf A}$ and ${\bf D}^2$ are centrosymmetric.

Theorem~\ref{theorem: H centrosymmetric} illustrates that Sinc basis functions preserve the parity property of the Sturm-Liouville operator when discretized. Hence, when the matrices $\,{\bf A}$ and ${\bf D}^2$ are symmetric centrosymmetric positive definite matrices, we can utilize these symmetries when solving for their generalized eigenvalues. In \cite{Cantoni1976a}, Cantoni et al. proved several properties of symmetric centrosymmetric matrices. In the following, we will utilize some of these properties to facilitate our task of obtaining approximations to the generalized eigenvalues of the matrices $\,{\bf A}$ and ${\bf D}^2$. The following lemma will demonstrate the internal block structure of symmetric centrosymmetric matrices.

\begin{lemma} \label{lemma: centrosymmetric, symmetric} \cite[Lemma 2]{Cantoni1976a}
If $\, {\bf H}$ is a square symmetric centrosymmetric matrix of dimension $\,(2N+1) \times (2N+1)$, then ${\bf H}$ can be written as:
\begin{equation}\label{formula: H matrix centrosymmet.ric}
{\bf H}  =  \left[\begin{array}{cccc} {\bf S} & {\bf x} & {\bf C^{T}}  \\
{\bf x^{T}} & h & {\bf x^{T}}{\bf J}\\
{\bf C} & {\bf J}{\bf x} & {\bf J}{\bf S}{\,\bf J}\end{array}\right],
\end{equation}
where ${\bf S}, {\bf C}$ are matrices of size $N \times N$, $ {\bf J}$ is the exchange matrix of size $N \times N$, ${\bf x} $ is a column vector of length $N$ and $h$ is a scalar. In addition, ${\bf S^{T}} = {\bf S}$ and  ${\bf C^{T}} = {\bf J}{\bf C}{\bf J}$.
\end{lemma}

The next lemma simplifies the calculation needed to solve for these eigenvalues.
\begin{lemma}\cite[Lemma 3]{Cantoni1976a}\label{lemma: Ortho similar}
Let ${\bf H}$ be a square symmetric centrosymmetric matrix as defined in lemma \ref{lemma: centrosymmetric, symmetric} and let ${\bf V}$ be a square matrix of dimension $(2N+1) \times (2N+1)$ defined by:
\begin{equation}
{\bf V}  = \left[\begin{array}{cccc} {\bf S-JC} & {\bf 0} & {\bf 0}  \\
{\bf 0} & h & \sqrt{2}~{\bf x^{T}}\\
{\bf 0} & \sqrt{2}~{\bf x} & {\bf S+JC}\end{array}\right] ,
\end{equation}
then ${\bf H} $ and ${\bf V}$ are orthogonally similar. That is, the matrices ${\bf H} $ and ${\bf V}$ have the same Jordan normal form and thus the same eigenvalue spectrum.
\end{lemma}

Cantoni et al. use Lemmas~\ref{lemma: centrosymmetric, symmetric} and \ref{lemma: Ortho similar} to prove the following Theorem concerning a standard eigenvalue problem where the matrix is centrosymmetric.

\begin{theorem} \label{theorem: symmetric centrosymmetric}\cite[Theorem 2]{Cantoni1976a}
Let ${\bf H}$ be a square symmetric centrosymmetric matrix as defined in Lemma \ref{lemma: centrosymmetric, symmetric}, then solving the eigenvalue problem $\det( {\bf H} - \lambda {\bf I} )=0$ is equivalent to solving the two smaller eigenvalue problems:
\begin{equation}
\det({\bf S-JC} -\lambda {\bf I} ) = 0  \quad \textrm{and} \quad
 \det \left( \left[ \begin{array}{cccc}  h & \sqrt{2}~{\bf x^{T}}\\
 \sqrt{2}~{\bf x} & {\bf S+JC}\end{array}\right] - \lambda \left[ \begin{array}{cccc}  1 & {\bf 0}\\
 {\bf 0} & {\bf I}\end{array}\right]  \right) = 0.
\end{equation}
\end{theorem}

Since our problem consists of solving a generalized eigenvalue problem where one matrix is a full symmetric centrosymmetric and the other is a diagonal centrosymmetric matrix, we propose the following Theorem.

\begin{theorem}\label{theorem: symmetric centrosymmetric generalized}
Let ${\bf H}$ and ${\bf W}$ be square symmetric centrosymmetric matrices of the same size, such that:
\begin{equation}
{\bf H}  =  \left[\begin{array}{cccc} {\bf S} & {\bf x} & {\bf C^{T}}  \\
{\bf x^{T}} & h & {\bf x^{T}}{\bf J}\\
{\bf C} & {\bf J}{\bf x} & {\bf J}{\bf S}{\,\bf J}\end{array}\right]  \qquad \textrm{and} \qquad
{\bf W}  =  \left[\begin{array}{cccc} \diag({\bf w}) & {\bf 0} & {\bf 0}  \\
{\bf 0} & w & {\bf 0}\\
{\bf 0} & {\bf 0} & {\bf J}\diag({\bf w}){\bf J}\end{array}\right] ,
\end{equation}
then solving the generalized eigenvalue problem $\det( {\bf H} - \lambda {\bf W} )=0$ is equivalent to solving the two smaller generalized eigenvalue problems:
\begin{equation}
\det({\bf S-JC} -\lambda \diag({\bf w}) ) = 0  \quad \textrm{and} \quad
 \det \left( \left[ \begin{array}{cccc}  h & \sqrt{2}~{\bf x^{T}}\\
 \sqrt{2}~{\bf x} & {\bf S+JC}\end{array}\right] - \lambda \left[ \begin{array}{cccc}  w & {\bf 0}\\
 {\bf 0} & \diag({\bf w})\end{array}\right]  \right) = 0.
\end{equation}
\end{theorem}

\underline{\bf Proof} \hskip 0.25cm This proof relies on the unitary transformation matrix presented in \cite[Lemma 3]{Cantoni1976a}:
\begin{equation}
{\bf K}  =  \dfrac{1}{\sqrt{2}}\left[\begin{array}{cccc} {\bf I} & {\bf 0} & -{\bf J}  \\
{\bf 0} & \sqrt{2} & {\bf 0}\\
{\bf I} & {\bf 0} & {\bf J}\end{array}\right],
\end{equation}
where ${\bf I}$ is the identity matrix and ${\bf J}$ is the exchange matrix.

It is easy to verify that:
\begin{equation}
{\bf K}{\bf H}{\bf K^{T}} = {\bf V},
\end{equation}
where ${\bf V}$ is the matrix in lemma \ref{lemma: Ortho similar}.

This result is analogous for the matrix ${\bf W}$ with a change in notation. Hence:
\begin{align}
0 = & \det( {\bf H} - \lambda {\bf W} ) \nonumber\\
 = & \det( {\bf K} ) \det( {\bf H} - \lambda {\bf W} )\det( {\bf K^{T}} ) \nonumber \\
  = &  \det( {\bf KHK^{T}} - \lambda {\bf KWK^{T}} ) \nonumber \\
  = & \det \left( \left[\begin{array}{cccc} {\bf S-JC} & {\bf 0} & {\bf 0} \nonumber \\
{\bf 0} & h & \sqrt{2}~{\bf x^{T}}\\
{\bf 0} & \sqrt{2}~{\bf x} & {\bf S+JC}\end{array}\right]- \lambda \left[\begin{array}{cccc} \diag({\bf w}) & {\bf 0} & {\bf 0}  \\
{\bf 0} & w & {\bf 0}\\
{\bf 0} & {\bf 0} & \diag({\bf w})\end{array}\right] \right) \nonumber \\
  = & \det({\bf S-JC} -\lambda \diag({\bf w}) )\det \left( \left[ \begin{array}{cccc}  h & \sqrt{2}~{\bf x^{T}}\\
 \sqrt{2}~{\bf x} & {\bf S+JC}\end{array}\right] - \lambda \left[ \begin{array}{cccc}  w & {\bf 0}\\
 {\bf 0} & \diag({\bf w})\end{array}\right]  \right),
\end{align}
from which the result follows.

Theorem \ref{theorem: symmetric centrosymmetric generalized} is very useful when $N$ is large since it is less costly to solve two symmetric generalized eigensystems of dimensions $N \times N$ and $(N+1) \times (N+1)$ rather than one symmetric eigensystem of dimension $(2N+1) \times (2N+1)$. Additionally, Lemma \ref{lemma: centrosymmetric, symmetric} also has large ramifications when it comes to saving storage space. As is discussed in \cite{Stenger1997a}, the $l^{\rm th}$ Sinc differentiation matrices are symmetric toeplitz matrices. Therefore, for a symmetric toeplitz matrix of dimension $(2N+1) \times (2N+1)$, only $2N+1$ elements need to be stored. Investigating the definition of the matrix ${\bf A}$ in equation \eqref{formula: A components}, we can see that ${\bf A}$ is defined as the sum of a symmetric toeplitz matrix and a diagonal matrix. Moreover, from Lemma \ref{lemma: centrosymmetric, symmetric} and Theorem \ref{theorem: symmetric centrosymmetric generalized}, using only the antidiagonal and anti-upper triangular half of matrix ${\bf C}$, the vector ${\bf x}$, the scalar $h$, the diagonal and lower triangular half of the matrix ${\bf S}$, the vector ${\bf w}$ and the scalar $w$, we can create all the elements needed to solve for the generalized eigenvalues of the matrices ${\bf A}$ and ${\bf D}^2$. Hence, the ratio of elements needed to be computed and stored at each iteration $N$ in order to solve for these eigenvalues is given by:
\begin{equation}
\textrm{Proportion of Entries Needed} \,=\, \dfrac{ \left(2N + N+1 \right) + \left(N +1\right) }{(2N+1)^2 + (2N+1)}  \,=\, \dfrac{1}{N+1}.
\end{equation}

Thus, only $\frac{1}{N+1}$  of the entries need to be generated and stored at every iteration to obtain all of the generalized eigenvalues.

In the following section, we will illustrate the gain in efficiency of these results by applying the DESCM to the Schr\"odinger equation with an anharmonic oscillator.

\section{The anharmonic oscillator}
The time independent Schr\"{o}dinger equation given by:
\begin{equation}\label{formula:Schrodinger equation}
{\cal H} \, \psi(x) \, \,=\,  E \, \psi(x) \qquad \textrm{with} \qquad  \lim_{|x|\to \infty}\psi(x) =0.
\end{equation}

In equation \eqref{formula:Schrodinger equation}, the Hamiltonian is given by the following linear operator:
$$
{\cal H} = -\dfrac{{\rm d}^2}{{\rm d} x^2} +V(x),
$$
where $V(x)$ is the potential energy function and $E$ is the energy eigenvalue of the hamiltonian operator ${\cal H}$. In our case, we are treating the anharmonic oscillator potential $V(x)$ defined by:
\begin{equation}\label{formula: anharmonic oscillator}
V(x)= \displaystyle \sum_{i=1}^{m}c_{i}x^{2i} \qquad  \textrm{with} \qquad  c_{m} > 0 \quad \textrm{and} \quad m \in \mathbb{N}\backslash\{1\}.
\end{equation}

In \cite{Safouhi49}, we successfully applied the DESCM to time independent Schr\"{o}dinger equation with an anharmonic potential. As we can see, the time independent Schr\"{o}dinger equation \eqref{formula:Schrodinger equation} is a special case of a Sturm-Liouville equation with $q(x)=V(x)$ and $\rho(x)=1$.  Applying Eggert et al.'s transformation and the DESCM with $\phi(x) =\sinh(x)$, we arrived at the following generalized eigenvalue problem:
\begin{equation}
\det({\bf A}-\mathcal{E}{\bf D}^{2}) = 0,
\label{formula: matrix solution-DET}
\end{equation}
where $\mathcal{E}$ are approximations of the energy eigenvalues $E$.

The matrices ${\bf A}$ and ${\bf D}^{2}$ defined by equation \eqref{formula: A components} and \eqref{formula: D components} are given by:
\begin{equation}\label{formula: A components Schrodinger}
A_{j,k} =   -\left(\dfrac{1}{h^{2}}\right)\delta^{(2)}_{j,k} \,  + \tilde{V}(kh)\delta^{(0)}_{j,k} \quad {\rm with} \quad -N \leq j,k \leq N,
\end{equation}
where:
\begin{equation}
\tilde{V}(x) \,=\, \dfrac{1}{4} - \dfrac{3}{4} \, \sech^{2}(x) + \cosh^{2}(x) \displaystyle \sum_{i=1}^{m}c_{i}\sinh^{2i}(x),
\end{equation}
and:
\begin{equation}\label{formula: D components Schrodinger}
D^2_{j,k} =   \cosh^{2}(kh) \, \delta^{(0)}_{j,k}  \qquad {\rm with} \qquad -N \leq j,k \leq N.
\end{equation}

Since the anharmonic potential $V(x)$ defined in \eqref{formula: anharmonic oscillator} is an even function, the function $\rho(x) =1$ is an even function and the conformal map $\phi(x) =\sinh(x)$ is an odd function, we know that Theorem \ref{theorem: H centrosymmetric} applies. Hence, the matrices ${\bf A}$ and ${\bf D}^{2}$ are symmetric centrosymmetric.

\section{Numerical Discussion}

In this section, we test the computational efficiency of the results obtained in Theorem \ref{theorem: symmetric centrosymmetric generalized}. All calculations are performed using the programming language Julia in double precision. The eigenvalue solvers in Julia use the linear algebra package {\it LAPACK}.

In \cite{Chaudhuri1991a}, Chaudhuri et al. presented several potentials with known analytic solutions for energy levels calculated using supersymmetric quantum mechanics, namely:
\begin{equation}\label{formula: true value energy}
\begin{array}{lllll}
V_{1}(x) & = & x^2 -4x^4+x^6  & \Rightarrow &  E_{0} = -2 \\
V_{2}(x) & = & 4x^2 -6x^4+x^6 & \Rightarrow &  E_{1} = -9 \\
V_{3}(x) & = & (105/64) x^2-(43/8)x^4 + x^6 -x^8 +x^{10} &\Rightarrow &  E_{0} = 3/8 \\
V_{4}(x) & = & (169/64)x^2 -(59/8)x^4 + x^6 -x^8 + x^{10}   &\Rightarrow &  E_{1} = 9/8.
\end{array}
\end{equation}

\begin{figure}[!ht]
\begin{center}
\begin{tabular}{cc} \includegraphics[width=0.35\textwidth]{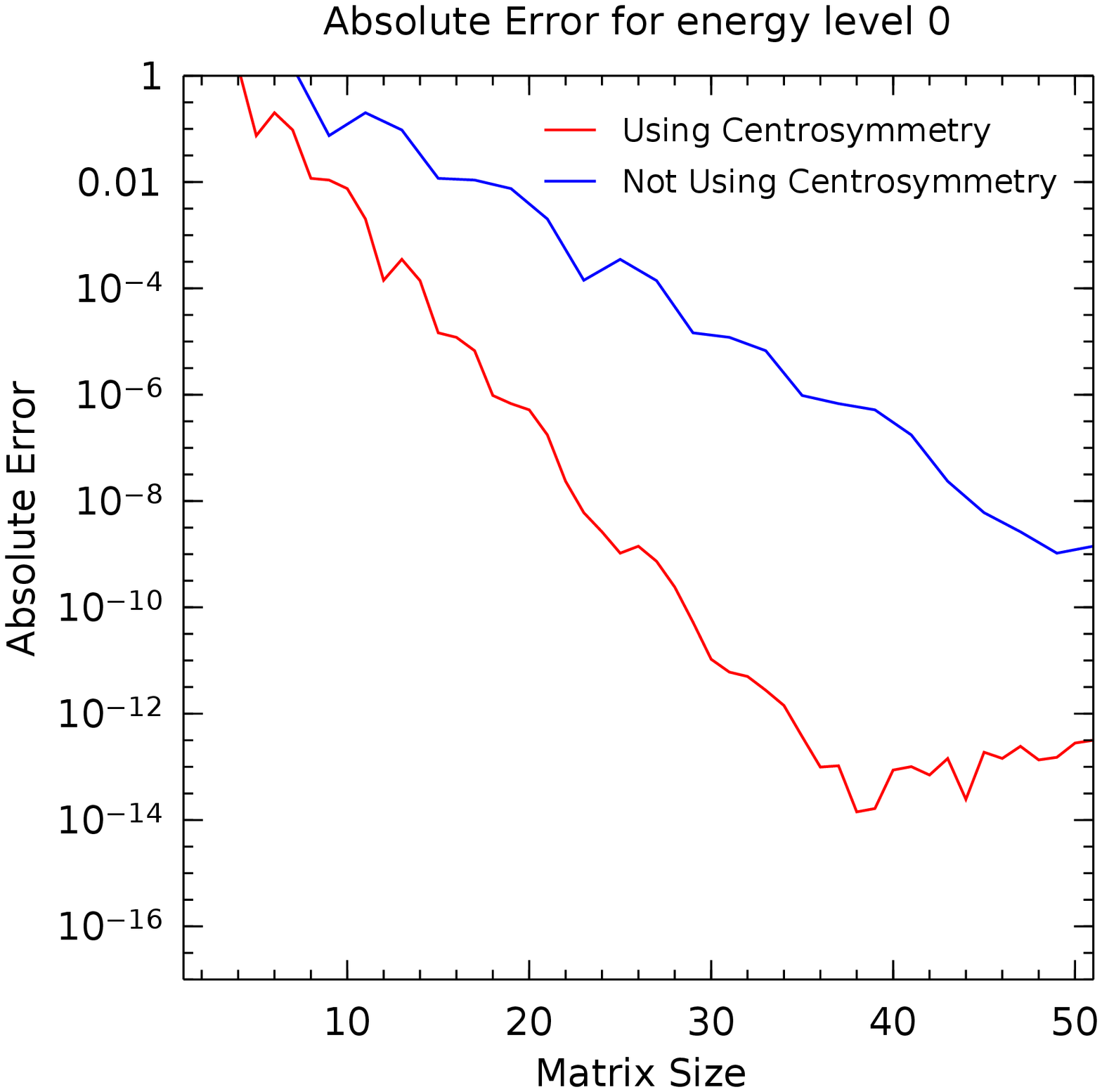} &  \includegraphics[width=0.35\textwidth]{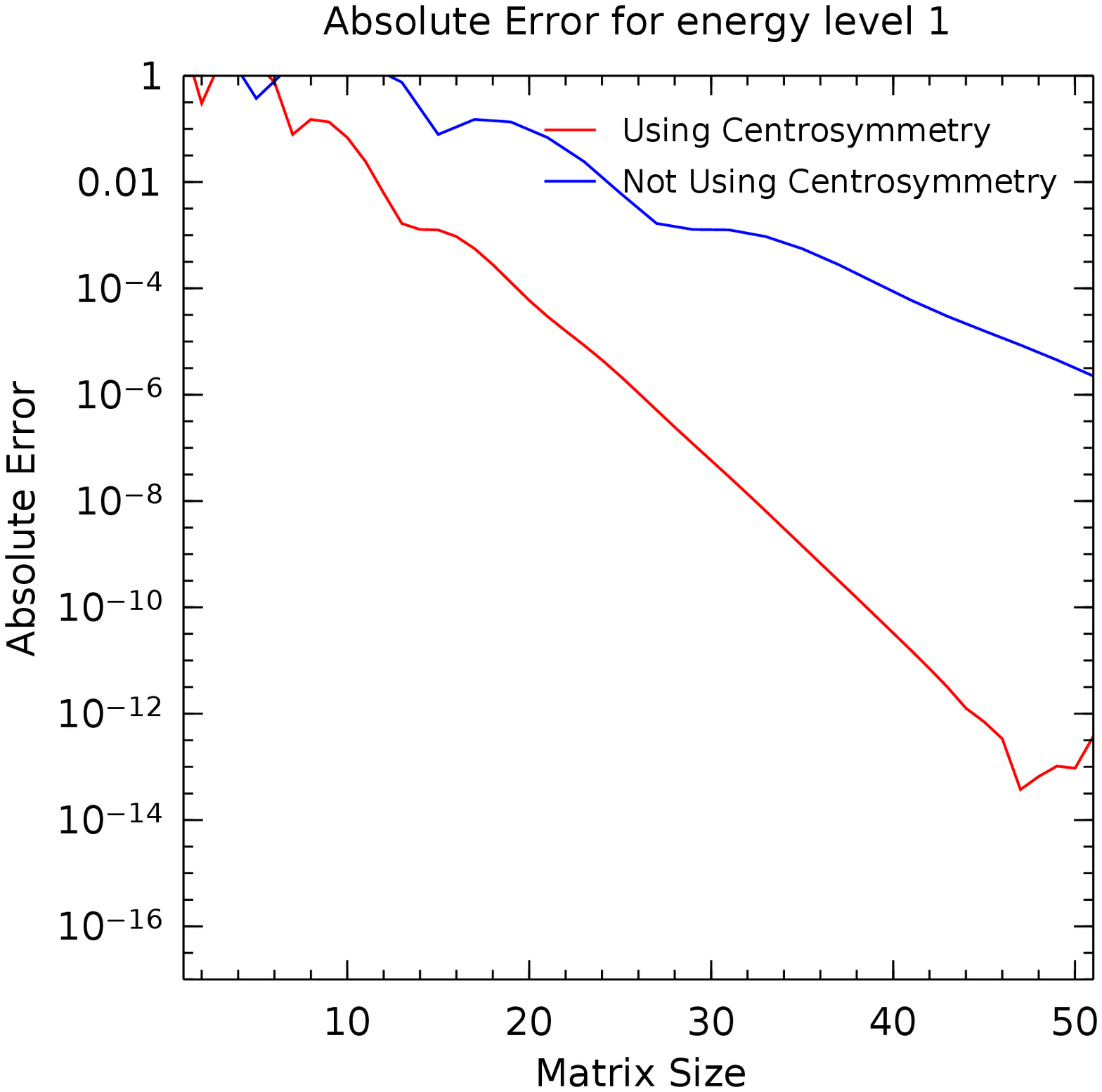} \\
(a) & (b) \\
\includegraphics[width=0.35\textwidth]{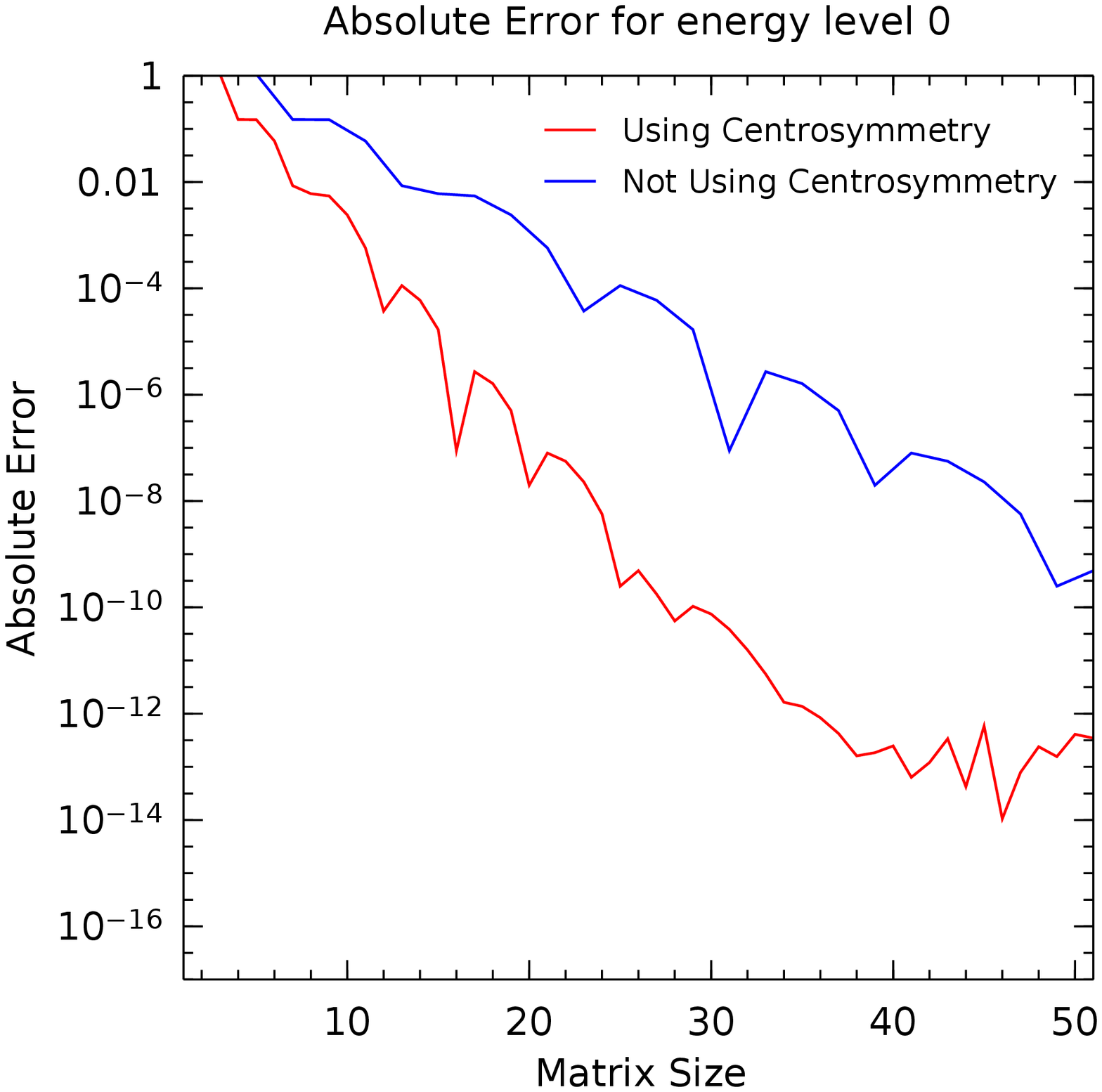} &  \includegraphics[width=0.35\textwidth]{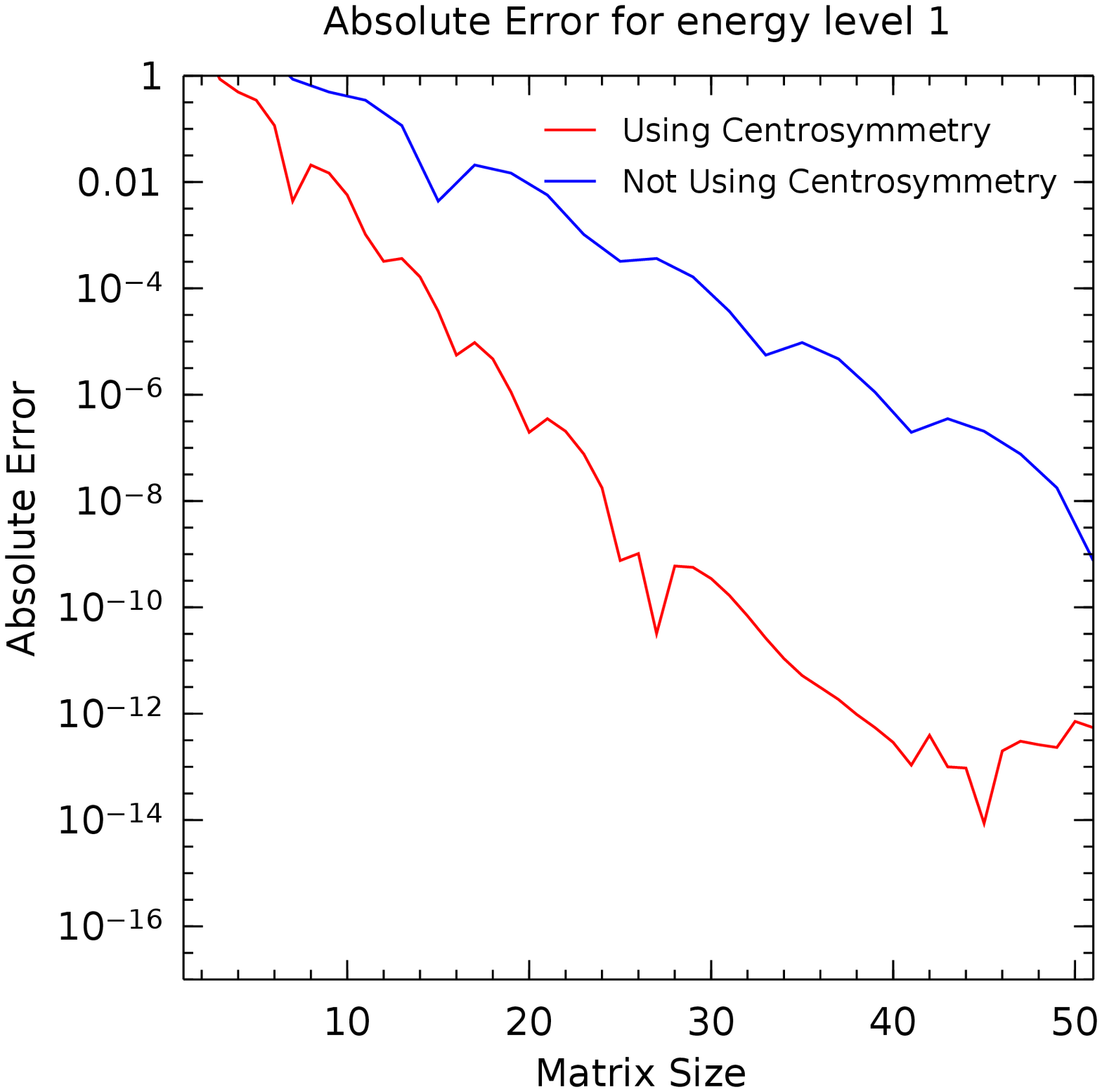} \\
(c) & (d)
\end{tabular}
\caption{Absolute error for the potentials $V_{i}(x)$ for $i=1,2,3,4$ given by~\eqref{formula: true value energy} with $\phi(x) = \sinh(x)$. \newline
(a)~$V_{1}(x) =  x^2 -4x^4+x^6$ with exact eigenvalue $E_{0} = -2$. (b)~$V_{2}(x)  =  4x^2 -6x^4+x^6$ with exact eigenvalue $E_{1} = -9$. (c)~$V_{3}(x)  =  (105/64) x^2-(43/8)x^4 + x^6 -x^8 +x^{10}$ with exact eigenvalue $E_{0} = 3/8$. (d)~$V_{4}(x)  =  (169/64)x^2 -(59/8)x^4 + x^6 -x^8 + x^{10} $ with exact eigenvalue $E_{1} = 9/8$.}
\label{figure: potentials}
\end{center}
\end{figure}

Figure~\ref{figure: potentials} presents the absolute error between our approximation and the exact values given in~\eqref{formula: true value energy}. The absolute error is defined by:
\begin{equation}
{\rm Absolute \, \, error} = \left| \mathcal{E}_{l}(N) - {\rm Exact \, \, value} \right| \qquad \textrm{for} \qquad l =0,1.
\end{equation}

The optimal mesh size obtained in \cite{Safouhi49}:
\begin{equation}\label{formula: Anharmonic optimal mesh size}
h = \dfrac{W \left( \frac{2^{m} \pi^{2}(m+1)N}{\sqrt{c_{m}}} \right)}{(m+1)N},
\end{equation}
where $W(x)$ is the Lambert W function, is used in the calculation.

As can be seen from Figure~\ref{figure: potentials}, using the centrosymmetric property improves the convergence rate of the DESCM significantly.

\section{Conclusion}

Sturm-Liouville eigenvalue problems are abundant in scientific and engineering problems. In certain applications, these problems possess a symmetry structure which results in the Sturm-Liouville operator to be commutative with the parity operator. As was proven in Theorem \ref{theorem: H centrosymmetric}, applying the DESCM  will preserve this symmetry and results in a generalized eigenvalue problem where the matrices are symmetric centrosymmetric. The centrosymmetric property leads to a substantial reduction in the computational cost when computing the eigenvalues by splitting the original eigenvalue problem of dimension $(2N+1)\times(2N+1)$ into two smaller generalized eigensystems of dimension $N\times N$ and $(N+1) \times (N+1)$. Moreover, due to the internal block structure of the matrices obtained using the DESCM, we have shown that only $\dfrac{1}{N+1}$ of all entries need to be computed and stored at every iteration in order to find all of their eigenvalues. Numerical results are presented for the time independent Schr\"odinger equation \eqref{formula:Schrodinger equation} with an anharmonic oscillator potential \eqref{formula: anharmonic oscillator}. Four exact potentials with known eigenvalues are tested and the results clearly demonstrated the reduction in complexity and increase in convergence.

\section{Tables and Figures}

\end{document}